\input amstex
\magnification=\magstep0
\documentstyle{amsppt}
\pagewidth{6.0in}
\input amstex
\topmatter
\title
Using the Jacobi-Trudi Formula to compute Stirling Determinants
\endtitle
\author
Tewodros Amdeberhan and Shalosh B. Ekhad \\
\endauthor
\subjclass 05A10, 05A15, 05E05,  11B65  \endsubjclass

\endtopmatter
\def\({\left(}
\def\){\right)}

\document

\noindent
\bf PART I: Theory \rm
\bigskip
\noindent
The \it unsigned Stirling numbers of the first kind \rm $\left[{n\atop k}\right]$ enumerate permutations of $n$ elements with $k$ disjoint cycles. They also arise as coefficients of the rising factorial, i.e.,
$$x(x+1)(x+2)\cdots(x+n-1)=\sum_{k=0}^n\left[{n\atop k}\right]x^k.$$
The \it Stirling numbers of the second kind \rm ${n \brace k}$ enumerate the number of ways to partition a set of $n$ objects into $k$ non-empty subsets. They also arise as coefficients of the falling factorial, i.e.,
$$\frac{x^k}{x(x-1)(x-2)\cdots(x-k)}=\sum_{n=1}^{\infty}{n \brace k}x^n.$$
Assume $a$ and $b$ to be non-negative integers. Our specific interest lies in computing the determinants of the following $n\times n$ matrices
$$M_n(a,b)=\left(\left[{i+a\atop j+b}\right]\right)_{1\leq i,j\leq n} \qquad \text{and} \qquad
N_n(a,b)=\left({i+a \brace j+b}\right)_{1\leq i,j\leq n}.$$ 
Denote $\beta_n(a,b)=\det(M_n(a,b))$ and $\gamma_n(a,b)=\det(N_n(a,b))$. Let $[a]=\{1,2,\dots,a\}$ be the integer interval. Given a partition $\lambda$, the \it Schur functions \rm can be given by
$$s_{\lambda}(\xi_1,\dots,\xi_a)=\frac{\det\left(\xi_i^{\lambda_j+a-j}\right)_{1\leq i,j\leq a}}{\det\left(\xi_i^{a-j}\right)_{1\leq i,j\leq a}}.$$
We are now ready to state our results.
\bigskip
\noindent
\bf Theorem 1. \it For $a, b, n \in\Bbb{Z}_{\geq0}$ and $b\leq a$, the sequence $\beta_n(a,b)$ has a rational generating function, in the variable $q$, with linearly factored denominator having the form
$$\prod\Sb i_1<i_2<\dots<i_b \\ i_1,i_2,\dots,i_b\in[a] \endSb\left(1-\frac{a!}{i_1i_2\cdots i_b}q\right).$$ \rm

\noindent
\bf Theorem 2. \it Denote $(n^c)=(n,n,\dots,n)\vdash cn$ to be a partition of $cn$. Then, we have
$$\beta_n(a,b)=s_{(n^{a-b})}(1,2,\dots,a) \qquad \text{and} \qquad \gamma_n(a,b)=s_{((a-b)^n}(1,2,\dots,a)).$$  \rm
\bf Proof. \rm First, observe that the Stirling numbers of the first and second kinds are the respective specializations of
$e_m$ and $h_m$. Both assertions follow from the Jacobi-Trudi and N\"agelsbach-Kostka identities
$$s_{\lambda}(\pmb{\xi})=\det(e_{\lambda_j'+j-i}(\pmb{\xi}))_{1\leq i,j\leq n}=\det(h_{\lambda_j+j-i}(\pmb{\xi}))_{1\leq i,j\leq n},$$
where $\lambda'$ is the \it conjugate \rm partition to $\lambda$ and $e_m(\pmb{\xi}), h_m(\pmb{\xi})$ are the elementary and the complete homogeneous symmetric functions, respectively.
$\square$

\bigskip
\noindent
\bf Theorem 3. \it Let $A_i=(-1)^{a-i}\cdot\frac{i^a}{a!}\binom{a}i$, for $i\in[a]$. Then, we have the "explicit" expressions
$$\align
\beta_n(a+b,b)=(-1)^{\binom{a}2}\sum\Sb 1\leq i_1<\cdots<i_a\leq a+b\endSb \,\,\prod_{\ell=1}^aA_{i_{\ell}}\cdot \prod \Sb \ell_u<\ell_v \\ \ell_u, \ell_v\in\{i_1,\dots,i_a\} \endSb 
(i_{\ell_u}-i_{\ell_v})^2\cdot \prod_{\ell=1}^a i_{\ell}^{n-a+1} =\gamma_a(n+b,b).
\endalign$$ \rm

\noindent
\bf Remark. \rm R. Stanley informed the first author that our identity for special case $\beta_n(a+1,a)$ looks like its should be equivalent to Exercise 7.4 in [2] (see also references therein).

\bigskip
\noindent
\bf PART II: Computations \rm
\bigskip
\noindent
Using the Jacobi-Trudi formula mentioned in Theorem 2 in Part I (Eq. (I.3.5) in [1], page 41) we computed explicit expressions for
$\beta_n(a,b)$, and also the explicit generating functions $\sum_{n=0}^{\infty} \beta_n(a,b) q^n$ (that are always rational functions of $q$),
for {\bf all} $10 \geq n \geq a \geq b \geq 0$.

\bigskip
\noindent
The output file is
\smallskip
{\tt https://sites.math.rutgers.edu/\~{}zeilberg/tokhniot/oStirlingDet1.txt}  

\bigskip
\noindent
It was generated by executing the command
\smallskip
{\tt Paper1(10,n,q):}
\smallskip
\noindent
in the Maple package accompanying this article, that can be gotten from
\smallskip

{\tt https://sites.math.rutgers.edu/\~{}zeilberg/tokhniot/StirlingDet.txt}

\bigskip

\hrule

\bigskip
\noindent
Tewodros Amdeberhan, Department of Mathematics, Tulane University, New Orleans, LA, 70118, USA. 
Email: {\tt tamdeber\@tulane.edu}

\bigskip
\noindent
Shalosh B. Ekhad (c/o D. Zeilberger) Department of Mathematics, Rutgers University \hfill\break
(New Brunswick), 
Hill Center-Busch Campus, 110 Frelinghuysen
Rd., Piscataway, NJ 08854-8019, USA. \hfill\break
Email: {\tt ShaloshBEkhad\@gmail.com}   
\bigskip
June 17, 2022
\bigskip

\enddocument